\documentclass[pamm,a4paper,fleqn%
]{w-art}
\usepackage{times,cite,w-thm}
\theoremstyle{plain}

\theoremstyle{definition}

\usepackage[]{graphicx}

\usepackage{amsfonts}
\usepackage{amssymb}

\def\Image{\hbox{\rm Image}}
\def\R{\mathbb R}

\def\dd{\hbox{d}}

\begin{document}
\DOIsuffix{theDOIsuffix}
\Volume{XX}
\Month{01}
\Year{2007}
\pagespan{1}{}



\title[On tangential stabilization in curvature driven flows of planar curves]
{On tangential stabilization in curvature driven flows of planar curves}


\author[K. Mikula]{Karol Mikula\inst{1}%
\address[\inst{1}]{Department of Mathematics, Slovak University of
Technology, Rad\-lin\-sk\'eho 11, 813 68 Bratislava, Slovak Republic}
}
\author[D. \v{S}ev\v{c}ovi\v{c}]{Daniel  \v{S}ev\v{c}ovi\v{c}\inst{2}
\address[\inst{2}]{Department of Applied Mathematics and Statistics,
FMFI, Comenius University, 842 48 Bratislava, Slovak Republic}
\footnote{Corresponding author: Daniel \v{S}ev\v{c}ovi\v{c}, \quad E-mail:~\textsf{sevcovic@fmph.uniba.sk},
            Phone: +00\,421\,2\,60295134
            Fax: +00\,421\,2\,65412305
\\
This research was supported by grants: VEGA 1/3321/06 (K.Mikula) 
and APVV-0247-06 (D.\v{S}ev\v{c}ovi\v{c}).
}
}

\begin{abstract}
We discuss the role of tangential stabilization in a curvature driven flow of 
planar curves. The governing system of nonlinear parabolic equations includes 
a nontrivial tangential velocity functional yielding a uniform redistribution 
of grid points along the evolving family of curves preventing numerically computed
curves from forming various instabilities. 
\end{abstract}
\maketitle                   





\section{Introduction}
In this paper we study evolution of a family of closed smooth plane curves
$\Gamma_t:S^1 \to \R^2$, $t\ge 0$, driven by the normal velocity $v$ which
is assumed to be a function of the curvature $k$, tangential angle $\nu$
and position vector $x\in\Gamma_t$,
\begin{equation}
v= \beta(x,k,\nu) \,.
\label{geomrov}
\end{equation}
As a typical example one can consider a normal velocity of the form: $v = k$ 
(mean curvature driven flow), $v = k^\frac13$ (affine invariant flow), 
$v= a(x,\nu) k + c(x,\nu)$  (Gibbs-Thomson law), etc. Geometric equations of the 
form (\ref{geomrov}) can often be found in variety of applied problems like e.g. 
the material science, dynamics of phase boundaries in thermomechanics, in modeling 
of flame front propagation, in combustion, in computations of first arrival times 
of seismic waves, in computational  geometry, robotics, semiconductors industry, etc. 
They also have a special conceptual importance in image processing and computer vision.
For an overview of important applications of (\ref{geomrov}) we refer to a
book by Sethian \cite{Se2}.

An idea behind the direct (or Langrangean) approach consists in representing the 
family of immersed curves $\Gamma_t$ by the position vector $x\in \R^2$, i.e.
$\Gamma_t=\Image(x(.,t))= \{x(u,t),\ u\in S^1\}$
where $x$ is a solution to the geometric equation
\begin{equation}
\partial_t x = \beta \vec N + \alpha \vec T
\label{pozrov}
\end{equation}
where $\beta=\beta(x,k,\nu)$, $\vec N=(-\sin \nu,\cos \nu)$
and $\vec T=(\cos \nu,\sin \nu)$ are the unit inward normal and tangent  vectors, 
respectively. We chose the orientation of the tangent vector $\vec T$ such that 
$\hbox{det}(\vec T,\vec N) =1$. Notice that the presence of arbitrary tangential
velocity functional $\alpha$ has no impact on the shape of evolving curves and thus 
$\alpha$ can be viewed as free parameter to be suitably determined. The unit 
arc-length parameterization of a curve $\Gamma=\Image(x)$ will be denoted by 
$s$. Then $d s = g\, d u$ where $g=|\partial_u x|$.

According to \cite{MS2,MS3} (see also \cite{MS4,MS5}) the system of governing equations 
for the curvature $k$, tangent angle $\nu$, local length $g$ and the position vector 
$x$ reads as follows: 
\begin{equation}
\partial_t k =\partial^2_s \beta +\alpha \partial_s k +k^2\beta\,,
\quad
\partial_t \nu = \beta^\prime_k \partial^2_s\nu
     + (\alpha + \beta^\prime_\nu) \partial_s\nu 
     + \nabla_x \beta.\vec T\,,
\quad
\partial_t g  = -g k\beta  + \partial_u\alpha\,,
\quad
\partial_t x  = \beta \vec N + \alpha \vec  T\qquad
\label{rovnice}
\end{equation}
where $(u,t)\in S^1\times(0,T)$, $\dd s=g\,\dd u$, A solution $(k,\nu,g,x)$ to (\ref{rovnice})  
is subject to initial conditions and periodic boundary conditions in the $u$ variable.

\section{The role of the tangential velocity functional}

Notice that the functional $\alpha$ is still undetermined and it may depend 
on  variables $k,\nu,g,x$ in various ways including  nonlocal dependence
in particular. Suitable choices of the tangential velocity functional $\alpha$ 
are discussed in a more detail in this section. Although $\alpha$ plays an 
important role in the governing equations resulting in dependence of $k,\nu,g,x$ 
on $\alpha$, the family of planar curves $\Gamma_t=\Image(x(.,t)), t\in[0,T),$ 
is independent of a particular choice of $\alpha$. 

\begin{figure}
\includegraphics[width=70mm]{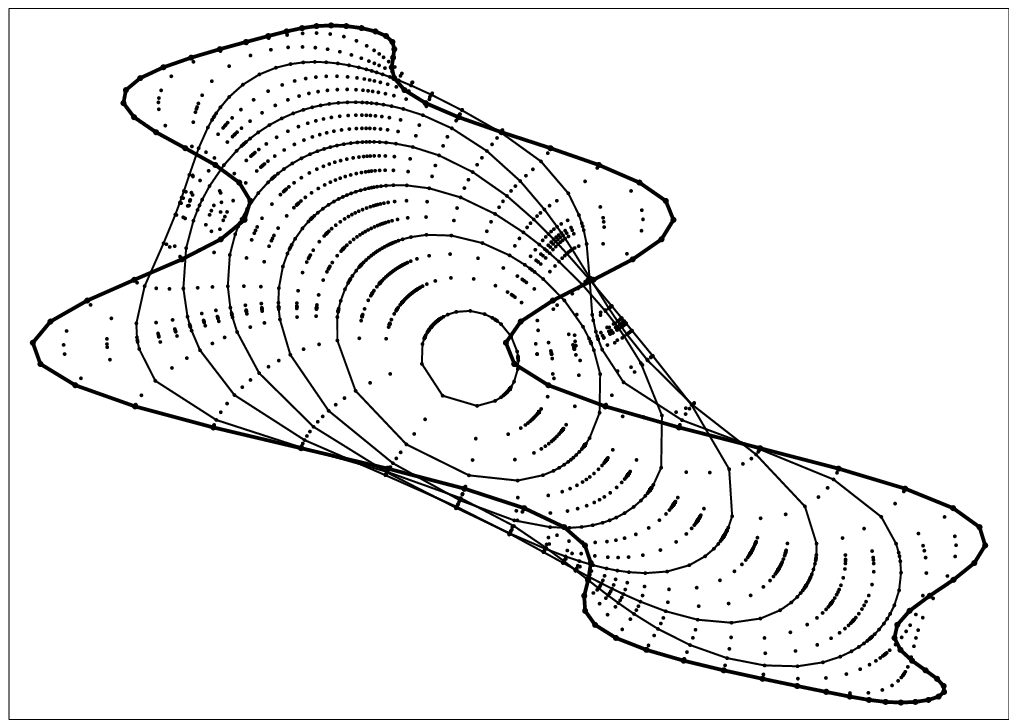}~a)
\hfil
\includegraphics[width=70mm]{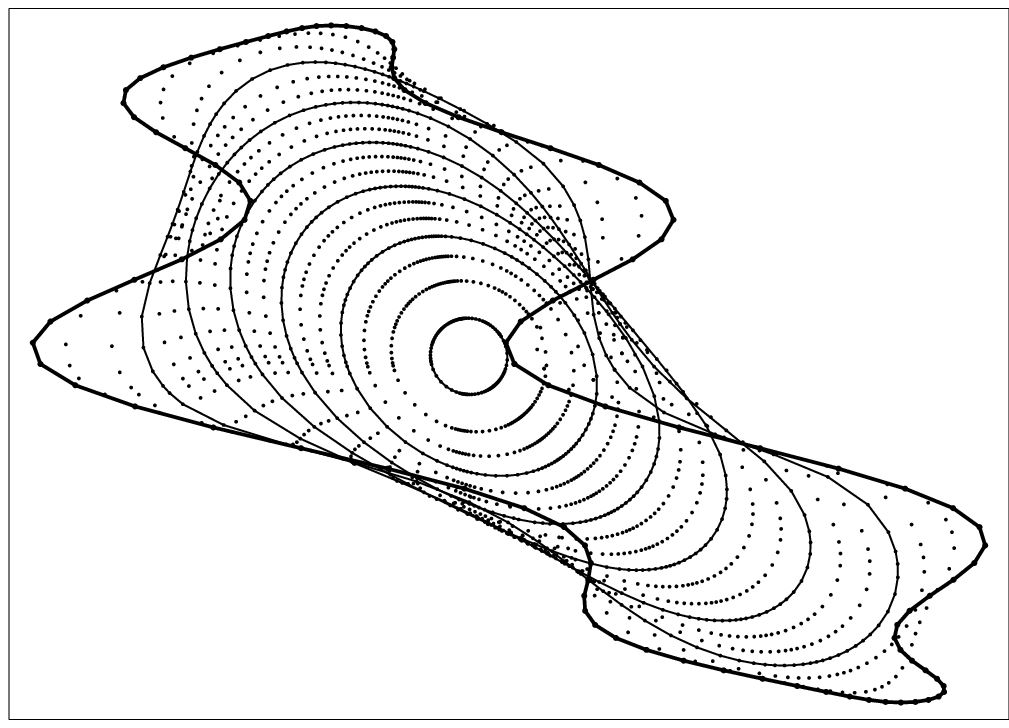}~b)
\caption{Numerically computed evolution of an initial curve (bold);  
a) merging of numerically computed grid points in the case of
zero tangential velocity $\alpha=0$; b) enhancement of grid point 
redistribution in the case of asymptotically uniform parametrization.}
\label{batman}
\end{figure}

To motivate further discussion, we  recall some of computational examples
in which the usual choice $\alpha=0$ fails and may lead to serious numerical
instabilities like e.g. formation of so-called swallow tails. In  Figure~\ref{batman}--a) 
we computed the mean curvature flow of an initial curve (bold faced curve). We chose 
$\alpha=0$. It should be obvious that numerically 
computed grid points merge in some parts of the curve $\Gamma_t$ preventing 
thus numerical approximation of $\Gamma_t, t\in[0,T),$ to be continued beyond
some time $T$ which is still far away from the maximal time of existence 
$T_{max}$. This and many other examples from \cite{MS2,MS3} showed that a 
suitable grid points redistribution  governed by a nontrivial tangential velocity 
functional $\alpha$ is needed in order to compute the solution over its life-span.

The idea behind construction of a suitable tangential velocity
functional $\alpha$ is rather simple and consists in the analysis of the
quantity $\theta$ defined as $\theta = \ln (g/L)$ where $g=|\partial_u x|$ 
is a local length and $L=L_t=\int_{\Gamma_t} \dd s=\int_0^1 g(u,t)\,\dd u$ is a total 
length of the curve  $\Gamma_t=\Image(x(.,t))$. The quantity $\theta$ can be 
viewed as the logarithm of the  relative local length ratio $g/L$. Taking into account 
equations (\ref{rovnice}) and the equation for the total length 
$\frac{d}{dt} L + \int_\Gamma k\beta \dd s =0$ (obtained again from  (\ref{rovnice}) 
by integration) we have
\begin{equation}
\partial_t \theta + k\beta -\langle k\beta \rangle_\Gamma 
= \partial_s \alpha\
\label{thetarovnica}
\end{equation}
where $\langle k \beta \rangle_{\Gamma}$ denotes the average  of 
$k \beta$ over the curve $\Gamma$, i.e. 
$\langle k \beta \rangle_{\Gamma} =  \frac{1}{L}\int_\Gamma k \beta \, \dd s$
By an appropriate choice of $\partial_s \alpha$ in the right hand side of 
(\ref{thetarovnica}) appropriately we can therefore control the behavior of 
$\theta$. Equation (\ref{thetarovnica}) can be also viewed as  a kind of a 
constitutive relation determining redistribution of grid point along a curve. 
The simplest possible choice of $\partial_s \alpha$ is:
\begin{equation}
\partial_s \alpha = k\beta -\langle k\beta \rangle_\Gamma
\label{prescribedredis}
\end{equation}
yielding $\partial_t\theta=0$ in (\ref{thetarovnica}). Consequently, 
$g(u,t)/L_t = g(u,0)/L_0 $ for any $u\in S^1,
\ t\in [0,T_{max})$. Notice that $\alpha$ can be uniquely computed from (\ref{prescribedredis})
under the additional renormalization constraint: $\alpha(0,t)=0$. The
tangential redistribution driven by a solution $\alpha$ to (\ref{prescribedredis}) 
is refereed to as {\it a parameterization preserving relative local length} 
(c.f. \cite{MS2}). It has been first discovered and utilized by Hou et al. 
in \cite{Hou1,Hou2} and independently by the authors in \cite{MS2,MS3}. 

A more general choice of $\alpha$ is based on the following setup:
\begin{equation}
\partial_s \alpha = k\beta -\langle k\beta \rangle_\Gamma
+ \left(  e^{-\theta} -1\right) \omega(t)
\label{uniformredis}
\end{equation}
where $\omega\in L^1_{loc}([0,T_{max}))$. If we additionally suppose
$\int_0^{T_{max}} \omega(\tau) \, \dd \tau = + \infty$
then, after insertion of  (\ref{uniformredis}) into (\ref{thetarovnica})  
and solving the ODE $\partial_t\theta =  \left(  e^{-\theta} 
-1\right) \omega(t),$ we obtain $\theta(u,t) \to 0$ as $t\to T_{max}$ and hence
$g(u,t)/L_t \to 1  \quad \hbox{as}\ t\to T_{max}$ uniformly w.r. to $u\in S^1$.
In this case redistribution of grid points along a curve becomes uniform as $t$ 
approaches the maximal time of existence $T_{max}$. We will refer to the 
parameterization based on (\ref{uniformredis}) to as {\it an asymptotically 
uniform parameterization} (c.f. \cite{MS3}). The impact of a tangential velocity functional 
defined as in (\ref{prescribedredis}) on enhancement of redistribution of grid 
points can be observed from two examples shown in Fig. \ref{batman}--b) computed 
by the authors in \cite{MS2}. It can be shown that the appropriate choice for the 
control function $\omega$ takes the form  $\omega = \kappa_1 + \kappa_2 \langle k\beta\rangle_\Gamma$ 
and $\kappa_1,\kappa_2\ge 0$ are given constants. A detailed discussion on this 
topic can be found in \cite{MS3,MS4}. If we insert tangential velocity functional 
$\alpha$ computed from (\ref{uniformredis}) into (\ref{rovnice}) the system of 
governing equations can be rewritten as follows:
\begin{equation}
\partial_t k =\partial^2_s \beta +\partial_s(\alpha k)  
+ k \langle k\beta \rangle_\Gamma  + \left(1-L/g\right) k \omega \,,
\qquad
\partial_t \nu = \beta^\prime_k \partial^2_s\nu
     + (\alpha + \beta^\prime_\nu) \partial_s\nu 
     + \nabla_x \beta.\vec T\,,
\label{prep-rovnice}
\end{equation}
\[
\partial_t g  = -g \langle k\beta   \rangle_\Gamma 
+ (L-g)\omega\,,
\qquad \partial_t x  = \beta \vec N + \alpha \vec T\,.  
\]
It is worth to note that the strong reaction term $k^2\beta$ in (\ref{rovnice}) has been 
replaced by the averaged term $k \langle k\beta \rangle_\Gamma$ in (\ref{prep-rovnice}). 
This is a very important feature as it allows for construction of an efficient and 
stable numerical scheme discussed in more details in \cite{MS3,MS4,MS5}.
 
%
\bibliographystyle{pamm}
\bibliography{sevcovic-mikula}
%

\end{document}